\documentclass{article}
\begin{document}
\newtheorem{proposition}{Proposition}[section]
\newtheorem{definition}{Definition}[section]
\newtheorem{lemma}{Lemma}[section]
\newcommand{\xl}{\stackrel{\rightharpoonup}{\cdot}}
\newcommand{\xr}{\stackrel{\leftharpoonup}{\cdot}}
\newcommand{\xlplus}{\stackrel{\rightharpoonup}{+}}
\newcommand{\xrplus}{\stackrel{\leftharpoonup}{+}}
\newcommand{\xluplus}{\stackrel{\rightharpoonup}{\uplus}}
\newcommand{\xruplus}{\stackrel{\leftharpoonup}{\uplus}}
\newcommand{\xlodot}{\stackrel{\rightharpoonup}{\odot}}
\newcommand{\xrodot}{\stackrel{\leftharpoonup}{\odot}}
\newcommand{\Ll}{\stackrel{\rightharpoonup}{L}}
\newcommand{\Lr}{\stackrel{\leftharpoonup}{L}}
\newcommand{\Rl}{\stackrel{\rightharpoonup}{R}}
\newcommand{\Rr}{\stackrel{\leftharpoonup}{R}}

\title{\bf Lying-Over Theorem \\on  Left Commutative Rngs}
\author{Keqin Liu\\Department of Mathematics\\The University of British Columbia\\Vancouver, BC\\
Canada, V6T 1Z2}
\date{November 11, 2005}
\maketitle

\begin{abstract} We introduce the notion of a graded integral element, prove the counterpart of the lying-over theorem on commutative algebra in the context of left commutative rngs, and use the Hu-Liu product to select a class of noncommutative rings.\end{abstract}

Left commutative rngs were introduced in \cite{Liu4}. I have two reasons to be interested in left commutative rngs. The first reason is that left commutative rngs are a class of commutative rings with zero divisors, and a graduate textbook about commutative algebra can be rewritten in the context of left commutative rngs if prime ideals are replaced by Hu-Liu prime ideals. Hence, studying left commutative rngs is an opportunity of extending the elegant commutative ring theory. Hu-Liu commutative trirings introduced in \cite{Liu5} are a class of noncommutative rings with zero divisors. Because of the first reason and the fact that 
left commutative rngs are special Hu-Liu commutative trirings, it seems to be true that Hu-Liu commutative trirings are suitable for reconsidering commutative ring theory. Hence, the new notions and ideas appearing in the study of left commutative rngs should be of much benefit to learning about the class of noncommutative rings with zero divisors. This is my second reason of being interested in left commutative rngs.

The purpose of this paper is to prove the lying-over theorem on left commutative rngs. My proof is based on the second proof of Theorem 3 on page 257 in \cite{ZS}. After avoiding the similar arguments, my proof is still much longer than the second proof in \cite{ZS}. This phenomenon is in fact hardly avoidable in the study of left commutative rngs. Therefore, it is predictable that the new textbook obtained by rewriting a textbook on commutative algebra in the context of left commutative rngs will be much thicker than the textbook on commutative algebra. 

After recalling basic definitions about left commutative rngs from \cite{Liu4}, we introduce the notion of a graded integral element in Section 1. Next, we prove the lying-over theorem on left commutative rngs in Section 2. Finally, we use the Hu-Liu product to select a class of noncommutative rings which have shown promise in the attempt to extend commutative ring theory.

Throughout this paper, we will use the notations in \cite{Liu4}. A ring means a ring with an identity, and a rng means a ring without an identity. 

\bigskip
\section{Basic Definitions}

Following \cite{Liu4}, we first define a left commutative rng in the following way.

\begin{definition}\label{def1.1} A rng $(\, R, \, + , \, \cdot \,)$ is called a {\bf left commutative rng} if $R$ satisfies the following four conditions.
\begin{description}
\item[(i)] The associative product $\cdot$ is {\bf left commutative}; that is, 
\begin{equation}\label{eq1}
xyz=yxz \quad\mbox{for $x$, $y$, $z\in R$.}
\end{equation}
\item[(ii)] There exists a {\bf left identity} $1^{\ell}$ of $R$ such that
\begin{equation}\label{eq2}
1^{\ell}x=x \quad\mbox{for $x\in R$.}
\end{equation}
\item[(iii)] There exists a binary operation $\,\sharp\,$ called the {\bf local product} on the {\bf additive halo} 
$$\hbar ^+(R): =\{\, x \,|\, \mbox{$x\in R$ and $x1^{\ell}=0$} \,\}$$ 
such that $(\, \hbar ^+(R) , \, + , \, \,\sharp\, \, )$ is a commutative ring with an identity $1^\sharp$. The identity $1^\sharp$ of the ring $\hbar ^+(R)$ is called the {\bf local identity} of $R$.
\item[(iv)] The two associative binary operations $\cdot$ and $\,\sharp\,$ satisfy the {\bf local strong Hu-Liu  triassociative law}:
\begin{equation}\label{eq3}
(x\alpha ) \,\sharp\, \beta=x(\alpha \,\sharp\, \beta),
\end{equation}
where $\alpha$, $\beta\in \hbar ^+(R)$ and $x\in R$.
\end{description}
\end{definition}

An equivalent way of defining a left commutative rng is the following

\begin{definition}\label{def1.2} Let $(\, R, \, + , \, \bullet_{_{1^{\ell}}} \,)$ be a commutative ring with an identity $1^{\ell}$. $R$ is called a {\bf left commutative rng} if there exists a binary operation $\cdot$ on $R$: $(x, \, y)\to xy:=x\cdot y$ such that 
$(\, R, \, + , \, \cdot \,)$ is a rng, and the conditions (ii), (iii), (iv) in Definition~\ref{def1.1} and the following hold:
\begin{equation}\label{eq4}
x\bullet_{_{1^{\ell}}}y=xy+yx-yx1^{\ell} \quad\mbox{for all $x$, $y\in R$.}
\end{equation}             
\end{definition}

The advantage of Definition~\ref{def1.2} is that it indicates explicitly that a left commutative rng is a commutative ring with many zero divisors. In fact, it follows from (\ref{eq4}) that 
$\alpha \bullet_{_{1^{\ell}}}\beta =0$ for all $\alpha$, $\beta\in \hbar ^+(R)$. Since 
$(\, R, \, + , \, \cdot \,)$ does not have an identity, $\hbar ^+(R)\ne 0$. Thus, every nonzero element of $\hbar ^+(R)$ is a zero divisor with respect to the commutative associative product 
$\bullet_{_{1^{\ell}}}$. The noncommutative associative product $\cdot$ is essential to the structure of a left commutative rng. Hence, we will use Definition~\ref{def1.1} in the study of left commutative rngs. 

\medskip
A left commutative rng $R$ is sometimes denoted by $(\, R, \, + , \, \cdot ,\, \,\sharp\, \,)$. Although $R$ may contain more than one left identity, we will fix a left identity $1^{\ell}$ of $R$ throughout this paper. Hence, it is meaningful to speak of the left identity of $R$.

\begin{definition}\label{def1.3} Let $(\, R, \, + , \, \cdot ,\, \,\sharp\, \,)$ be a left commutative rng with the left identity $1^{\ell}$, and let $I$ be a subgroup of 
$(\, R, \, +  \,)$.
\begin{description}
\item[(i)] $I$ is called an {\bf ideal} of $R$ if $IR\subseteq I$, $RI\subseteq I$ and 
$I\cap \hbar ^+(R)$ is an ideal of the ring $(\, \hbar ^+(R) , \, + , \, \,\sharp\, \, )$.
\item[(ii)] $I$ is called a {\bf left commutative subrng} of $R$ if $1^{\ell}\in I$, 
$II\subseteq I$, and $I\cap \hbar ^+(R)$ is subring of the ring $(\, \hbar ^+(R) , \, + , \, \,\sharp\, \, )$.
\end{description}
\end{definition}

Let $(\, R, \, + , \, \cdot ,\, \,\sharp\, \,)$ be a left commutative rng with the left identity $1^{\ell}$. Then we have
\begin{equation}\label{eq5}
R=R_0 \oplus R_1 \quad\mbox{(as Abelian groups)},
\end{equation}
where $R_0 : =R1^{\ell}$ and $R_1 : = \hbar ^+(R)$. $R_{\varepsilon}$ is called the 
{\bf $\varepsilon$-part} of $R$ for $\varepsilon=0$ and $1$. Every element $A$ of $R$ can be expressed uniquely as
\begin{equation}\label{eq6}
a=a_0 + a_1 \quad\mbox{with $a_0:=a1^{\ell}$ and $a_1:=a-a1^{\ell}$},
\end{equation}
and $a_{\varepsilon}$ is called the {\bf $\varepsilon$-component} of $a$, where $\varepsilon=0$ and $1$. If $I$ is an ideal or a left commutative subrng of $R$, then $I$ respects (\ref{eq5});
that is,
\begin{equation}\label{eq7}
I=I_0 \oplus I_1 \quad\mbox{with $I_{\varepsilon}:=I\cap R_{\varepsilon}$.}
\end{equation}

\medskip
The next definition is based on Proposition 1.1 (i) in \cite{Liu4}.

\begin{definition}\label{def1.4} Let $(\, R, \, + , \, \cdot ,\, \,\sharp\, \,)$ be a left commutative rng with the left identity $1^{\ell}$, and let $p=p_0\oplus p_1\ne R$ be an ideal of $R$, where $p_{\varepsilon}:=p\cap R_{\varepsilon}$, $R_0 : =R1^{\ell}$ and $R_1 : = \hbar ^+(R)$. $p$ is called a {\bf Hu-Liu prime ideal} of $R$ if
\begin{equation}\label{eq8}
\mbox{$x_0y_{\varepsilon}\in p  \,\Rightarrow \, $ $x_0\in p_0$ or $y_{\varepsilon}\in p_{\varepsilon}$}
\end{equation}
and
\begin{equation}\label{eq9}
\mbox{$x_1\,\sharp\, y_1\in P_1\,\Rightarrow \,$  $x_1\in P_1$ or $y_1\in P_1$,}
\end{equation}
where $x_{\varepsilon}$, $y_{\varepsilon}\in R_{\varepsilon}$ and $\varepsilon=0$, $1$.
\end{definition}

The set of all Hu-Liu prime ideals of $R$ is called the {\bf spectrum} of $R$ and is denoted by $spec^{\sharp}R$.

\medskip
We finish this section by introducing the notion of a graded integral element.

\begin{definition}\label{def1.5} Let $R$ be a left commutative subrng of a left commutative rng $U$. \begin{description}
\item[(i)] An element $u$ of $U$ is said to be {\bf graded integral over $R$} if $u_{\varepsilon}$ is integral over the ring $R_{\varepsilon}$ for $\varepsilon=0$ and $1$, where $u_{\varepsilon}$ is the $\varepsilon$-component of $u$, $R_0 : =R1^{\ell}$, $R_1 : = \hbar ^+(R)$ and $1^{\ell}$ is the left identity of $U$.
\item[(ii)] $U$ is said to be {\bf graded integral over $R$} if every element $u$ of $U$ is graded integral over $R$.
\end{description}
\end{definition}

\bigskip
\section{The Lying-Over Theorem on Left Commutative Rngs}

This section devotes to the proof of the following

\begin{proposition}\label{pr1.1} ({\bf The Lying-Over Theorem on Left Commutative Rngs}) Let $R$ be a left commutative subrng of a left commutative rng $U$, and suppose that $U$ is graded integral over $R$. If $p\in spec^{\sharp}R$, then there exists $q\in spec^{\sharp}U$ such that $q\cap R=p$.
\end{proposition}

\medskip
\noindent
{\bf Proof} Consider the set
$$
T:=\{\, J \,|\, \mbox{$J$ is an ideal of $U$ and $J\cap R\subseteq p$}\,\}.
$$

Since the partial order set $(\, T, \, \subseteq \,)$ satisfies the hypotheses of Zorn's Lemma, $T$ has a maximal element $q$. We are going to prove 
\begin{equation}\label{eq10} 
q\cap R=p
\end{equation}
and
\begin{equation}\label{eq11} 
\mbox{$q$ is a Hu-Liu prime ideal of $U$}.
\end{equation}

\medskip
Suppose that $q\cap R\ne p$. Then $q\cap R\subset p$. Since $q\cap R$ is an ideal of $R$, 
$q\cap R$ respects (\ref{eq5}). Hence, we have
\begin{equation}\label{eq12} 
(q\cap R)_0\oplus (q\cap R)_1=q\cap R\subset p=p_0\oplus p_1,
\end{equation}
where $p_{\varepsilon}=p\cap R_{\varepsilon}$ and
$$
(q\cap R)_{\varepsilon}=(q\cap R)\cap R_{\varepsilon}=(q\cap U_{\varepsilon})\cap R_{\varepsilon}
=q_{\varepsilon}\cap R_{\varepsilon}.
$$

It follows from (\ref{eq12}) that either $q_0\cap R_0\subset p_0$ or $q_1\cap R_1\subset p_1$.

\medskip
\underline{\it Case 1:} $q_0\cap R_0\subset p_0$, in which case, there exists $x_0\in R_0$ such that
\begin{equation}\label{eq13} 
x_0\in p_0 \quad\mbox{and}\quad x_0\not\in q_0.
\end{equation}
Hence, $q\subset q+x_0U$. Since $q+x_0U$ is an ideal of $U$, we have 
$(q+x_0U)\cap R\not\subseteq p$, which implies that
\begin{equation}\label{eq14} 
(q_{\varepsilon}+x_0U_{\varepsilon})\cap R_{\varepsilon}=
\Big((q+x_0U)\cap R\Big)_{\varepsilon}\not\subseteq p_{\varepsilon} 
\quad\mbox{for $\varepsilon=0$ or $1$}.
\end{equation}

If (\ref{eq14}) is true for $\varepsilon=1$, then there exists $a_1\in q_1$, $u_1\in U_1$ and 
$r_1\in R_1\setminus p_1$ such that $a_1+x_0u_1=r_1$ or
\begin{equation}\label{eq15} 
x_0u_1\equiv r_1 (mod\,\, q_1).
\end{equation}

Since the ring $(\, U_1, \, + , \, \,\sharp\, \,)$ is integral over the ring 
$(\, R_1, \, + , \, \,\sharp\, \,)$, we have
\begin{equation}\label{eq16} 
u_1^{\,\sharp\, n}+a_{11}\,\sharp\, u_1^{\,\sharp\, (n-1)}+\cdots + a_{n-1, 1}\,\sharp\, u_1
+a_{n1}=0
\end{equation}
for some $n\in\mathcal{Z}_{>0}$, $a_{11}$, $\dots$, $a_{n1}\in R_1$. Multiplying (\ref{eq16}) from the left side by $x_0^n$ and using (\ref{eq3}), we get
\begin{equation}\label{eq17} 
(x_0u_1)^{\,\sharp\, n}+(x_0a_{11})\,\sharp\, (x_0u_1)^{\,\sharp\, (n-1)}+\cdots + 
(x_0^{n-1}a_{n-1, 1})\,\sharp\, (x_0u_1)
+x_0^na_{n1}=0.
\end{equation}

By (\ref{eq15}), (\ref{eq17}) gives the following
\begin{equation}\label{eq18} 
r_1^{\,\sharp\, n}+(x_0a_{11})\,\sharp\, r_1^{\,\sharp\, (n-1)}+\cdots + 
(x_0^{n-1}a_{n-1, 1})\,\sharp\, r_1
+x_0^na_{n1}\equiv 0 (mod\,\, q_1).
\end{equation}
Note that the left side of (\ref{eq18}) is in $R_1$ and $q_1\cap R_1\subseteq p_1$. Thus, (\ref{eq18}) can be rewritten as
\begin{equation}\label{eq19} 
r_1^{\,\sharp\, n}+(x_0a_{11})\,\sharp\, r_1^{\,\sharp\, (n-1)}+\cdots + 
(x_0^{n-1}a_{n-1, 1})\,\sharp\, r_1
+x_0^na_{n1}\equiv 0 (mod\,\, p_1).
\end{equation}

Using (\ref{eq13}) and the fact that $p$ is a Hu-Liu prime ideal of $R$, we get from (\ref{eq19}) that $r_1^{\,\sharp\, n}\in p_1$ or $r_1\in p_1$, which is impossible.

If (\ref{eq14}) is true for $\varepsilon=0$, then we also get a contradiction by the similar argument above.

\medskip
\underline{\it Case 2:} $q_1\cap R_1\subset p_1$, in which case, a contradiction is derive from using the ordinary Lying-Over Theorem on commutative algebra or using the same argument as the one in Case 1.

\medskip
This proves (\ref{eq10}). We now begin to prove (\ref{eq11}).

First, if $x_0$, $y_0\in U_0\setminus q_0$, then
\begin{equation}\label{eq20} 
q\subset q+x_0U \quad\mbox{and}\quad q\subset q+y_0U
\end{equation}
Since $p=q\cap R\subseteq (q+x_0U)\cap R$ and $p=q\cap R\subseteq (q+y_0U)\cap R$, (\ref{eq20}) implies that
\begin{equation}\label{eq21} 
p\subset (q+x_0U)\cap R \quad\mbox{and}\quad p\subset (q+y_0U)\cap R,
\end{equation}
which imply that there exist $a_{\varepsilon}$, $b_{\varepsilon}\in q_{\varepsilon}$,
$u_{\varepsilon}$, $v_{\varepsilon}\in U_{\varepsilon}$ and
$s_{\varepsilon}$, $t_{\varepsilon}\in R_{\varepsilon}\setminus p_{\varepsilon}$ such that
\begin{equation}\label{eq22} 
a_{\varepsilon}+x_0u_{\varepsilon}=s_{\varepsilon} \quad\mbox{and}\quad 
b_{\delta}+y_0u_{\delta}=t_{\delta}
\end{equation}
for some $\varepsilon$, $\delta\in \{0,\, 1\}$. We are going to prove
\begin{equation}\label{eq23} 
x_0y_0\not\in q_0.
\end{equation}

By (\ref{eq22}), we need to consider four cases. The proofs in the four cases are similar. As an example, let us consider the case: (\ref{eq22}) holds for $\varepsilon=1$ and $\delta=0$. In this case, we have $a_1+x_0u_1=s_1$ and $b_0+y_0v_0=t_0$. It follows that
\begin{eqnarray*}
x_0y_0\in q_0&\Rightarrow& t_0s_1=(b_0+y_0v_0)(a_1+x_0u_1)\\
&=&\underbrace{b_0a_1+b_0x_0u_1+y_0v_0a_1}_{\mbox{This is an element of $q_1$}}+
(x_0y_0)v_0u_1\in q_1\\
&\Rightarrow& t_0s_1\in q_1\cap R_1=p_1\subseteq p\\
&\Rightarrow& t_0\in p_0 \quad\mbox{or}\quad s_1\in p_1,
\end{eqnarray*}
which contradicts the choice of $t_0$ or $s_1$. Hence, (\ref{eq23}) holds. This proves that
\begin{equation}\label{eq24} 
x_0, \, \, y_0\in U_0\setminus q_0 \Rightarrow x_0y_0\in U_0\setminus q_0.
\end{equation}

\medskip
Next, if $x_0\in U_0\setminus q_0$ and $y_1\in U_1\setminus q_1$, then
$$
q\subset q+x_0U \quad\mbox{and}\quad q\subset q+y_1\,\sharp\, U_1,
$$
which imply that
$$
p\subset(q+x_0U)\cap R  \quad\mbox{and}\quad p\subset(q+y_1\,\sharp\, U_1)\cap R.
$$
Hence, there exist $a_{\varepsilon}\in q_{\varepsilon}$, $u_{\varepsilon}\in U_{\varepsilon}$,
$s_{\varepsilon}\in R_{\varepsilon}\setminus p_{\varepsilon}$, $c_1\in q_1$, $w_1\in U_1$ and
$r_1\in R_1\setminus p_1$ such that
\begin{equation}\label{eq25} 
a_{\varepsilon}+x_0u_{\varepsilon}=s_{\varepsilon} \quad\mbox{and}\quad 
c_1+y_1\,\sharp\, w_1=r_1
\end{equation}
for some $\varepsilon =0$ or $1$.

If $a_0+x_0u_0=s_0$ and $c_1+y_1\,\sharp\, w_1=r_1$, then
\begin{eqnarray*}
x_0y_1\in q_1&\Rightarrow& s_0r_1=(a_0+x_0u_0)(c_1+y_1\,\sharp\,w_1)\\
&=&\underbrace{a_0c_1+a_0(y_1\,\sharp\,w_1)+x_0u_0c_1}_{\mbox{This is an element of $q_1$}}+
\Big(u_0(x_0y_1)\Big)\,\sharp\,w_1\in q_1\\
&\Rightarrow& s_0r_1\in q_1\cap R_1=p_1\subseteq p\\
&\Rightarrow& s_0\in p_0 \quad\mbox{or}\quad r_1\in p_1,
\end{eqnarray*}
which is impossible. This proves that
\begin{equation}\label{eq26} 
x_0\in U_0\setminus q_0 \quad\mbox{and}\quad y_1\in U_1\setminus q_1\Rightarrow 
x_0y_1\in U_1\setminus q_1.
\end{equation}

If (\ref{eq25}) is true for $\varepsilon=1$, then (\ref{eq26}) can be proved similarly.

\medskip
Finally, we have
\begin{equation}\label{eq27} 
x_1, \,\, y_1\in U_1\setminus q_1\Rightarrow x_1\,\sharp\,y_1\in U_1\setminus q_1.
\end{equation}

\medskip
It follows from (\ref{eq24}), (\ref{eq26}) and (\ref{eq27}) that (\ref{eq11}) holds.

\medskip
This completes the proof of Proposition~\ref{pr1.1}.

\hfill\raisebox{1mm}{\framebox[2mm]{}}

\bigskip
\section{Hu-Liu Commutative Rings}

First, we introduce the following

\begin{definition}\label{def3.1} Let $(\, R, \, + , \, \bullet \,)$ be a ring with an identity $\sigma$. The ring product $\bullet$ is called a {\bf Hu-Liu product} induced by $\sigma$, and the ring $R$ is called a {\bf ring with the Hu-Liu product} if the following three conditions are satisfied:
\begin{description}
\item[(i)] There exist two binary operations $\xl$ and $\xr$ such that
\begin{equation}\label{eq1.201}
x\bullet _{\sigma}y=x\xl y +x\xr y- (x\xr \sigma)\xl y
\end{equation}
for all $x$, $y\in R$. 
\item[(ii)] The three operations $\bullet _{\sigma}$, $\xl$ and $\xr$ satisfy the {\bf strong Hu-Liu triassocitive law}:
\begin{eqnarray}
\label{eq1.81} (x\xl y)\bullet_{\sigma} z&=&x\bullet_{\sigma} (y\xr z),\\
\label{eq1.84} x\xl (y\bullet_{\sigma} z)&=&(x\xl y)\xl z,\\
\label{eq1.85} (x\bullet_{\sigma} y)\xr z&=&(x\xr y)\xr z
\end{eqnarray}
for all $x,y,z\in R$.
\item[(iii)] The distributive laws
$$ x\ast (y+z)=x\ast y +x\ast z, \quad (y +z)\ast x= y\ast x + z\ast x $$
hold for all $x$, $y$, $z\in R$ and $\ast\in \{\, \xl ,\, \xr \, \}$.
\end{description}
\end{definition}

Note that we do not assume that any of the two binary operations $\xl$ and $\xr$ is associative in Definition~\ref{def3.1}. However, by Proposition 1.2 and 1.3 in \cite{Liu5}, each of the two binary operations $\xl$ and $\xr$ must be associative, and  the two binary operations $\xl$ and $\xr$ satisfy the diassociative law.

We will use $(\, R, \, + , \, \bullet \,)$ or $(\, R, \, + , \, \bullet ,\, \xl ,\, \xr \,)$ to denote a ring with the Hu-Liu product $\bullet$. 

\begin{definition}\label{def3.2} Let $(\, R, \, + , \, \bullet ,\, \xl ,\, \xr\,)$ be a ring with the Hu-Liu product. $R$ is called a {\bf Hu-Liu commutative ring} if
$$
x\xl y -y\xr x \,\in\,\hbar ^+ (R) \quad\mbox{for all $x$, $y\in R$},
$$
where $\hbar ^+ (R)=\{\, x \in R \, | \, \sigma\xl x=0 \,\}$ and $\sigma$ is the identity of the ring $R$.
\end{definition}

The field-like objects appearing in my research come from Hu-Liu commutative ring.

\end{document}